\documentclass{amsart}
\usepackage{amssymb,latexsym}
\usepackage{amsmath}
\usepackage[dvipdfm]{graphicx}
\theoremstyle{plain}
\newtheorem{thm}{\bfseries Theorem}
\newtheorem{lemma}[thm]{\bfseries Lemma}        
\newtheorem{remark}[thm]{\bfseries Remark}    

\newtheorem{fact}[thm]{\bfseries Fact}

\newtheorem{prob}[thm]{\bfseries Problem}

\begin{document}

\title[Elementary proof techniques for islands]{Elementary proof techniques for the maximum number of islands}

\author[J. Bar\'at]{J\'anos Bar\'at}
\address{Bolyai Institute, University of Szeged, Aradi v\'ertan\'uk tere 1, 6720 Szeged, Hungary}
\curraddr{Department of Computer Science, University of Pannonia, Egyetem u. 10, 8200 Veszpr\'em, Hungary}
\thanks{Research is supported by OTKA Grant PD~75837.}

\author[P. Hajnal]{P\'eter Hajnal}
\address{Bolyai Institute, University of Szeged, Aradi v\'ertan\'uk tere 1, 6720 Szeged, Hungary}
\thanks{Research is supported by OTKA Grant K~76099.}

\author[E. K. Horv\'ath]{Eszter K. Horv\'ath}
\address{Bolyai Institute, University of Szeged, Aradi v\'ertan\'uk tere 1, 6720 Szeged, Hungary}
\thanks{This research was partially supported by OTKA Grant T~049433, and by the Provincial Secretariat
for Science and Technological Development, Autonomous  Province of
Vojvodina, grant "Lattice methods and applications". }




\date{\today}

\begin{abstract}
Islands are combinatorial objects that can be intuitively defined on a board
consisting of a finite number of cells.
Based on the neighbor relation of the cells, it is a fundamental property
that two islands are either containing or disjoint.
Recently, numerous extremal questions have been answered using different methods.
We show elementary techniques unifying these approaches.
Our building parts are based on rooted binary trees and discrete geometry.

Among other things, we show the maximum cardinality of islands on a toroidal board
and in a hypercube.
We also strengthen a previous result by rarefying the neighborhood relation.
\end{abstract}

\maketitle

\section{Introduction, preliminaries}\label{intro}

We start with an intuitive notion.
Let a rectangular $m\times n$ board be given.
We associate a number (real or integer) to each cell of the board.
We can think of this number as a height above see level.
A rectangular part of the board is called a  {\it rectangular island},
if and only if there is a possible water level such that the rectangle
is an island in the usual sense.

\begin{figure}[ht]
 \begin{center}
  \includegraphics[scale=0.5]{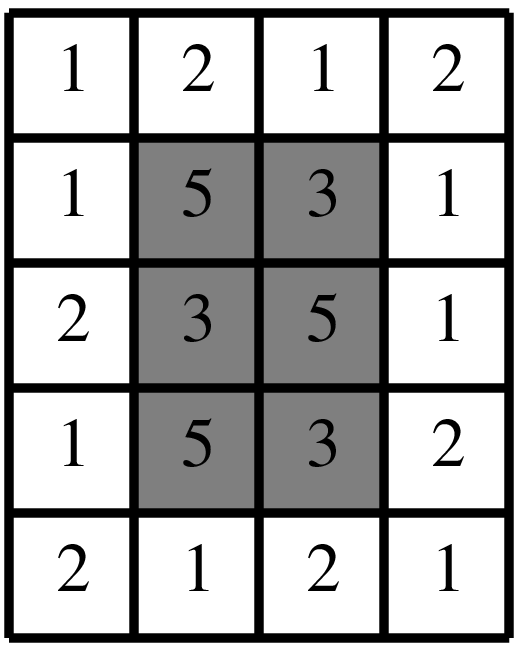}
 \end{center}
   \caption{Rectangular landscape with heights}
  \label{island}
\end{figure}

The notion of an island turned up recently in information theory.
The characterization of the lexicographical length sequences of binary maximal instantaneous
codes in \cite{FS} uses the notion of {\it full segments}, which are one-dimensional islands.
Several generalizations led to interesting combinatorial problems.
G. Cz\'edli discovered a connection between islands and weakly independent subsets of finite distributive lattices.
He determined the maximum number of rectangular islands on a rectangular board \cite{Cz}.
Cz\'edli's method is based on weak bases of a finite distributive lattice \cite{CzHSch}.
G. Pluh\'ar \cite{P} gave upper and lower bounds in higher dimensions.
E. K. Horv\'ath, Z. N\'emeth and G. Pluh\'ar \cite{HNP} gave upper and lower bounds for the maximum number of triangular islands
on a triangular board.
In \cite{L} the minimal size of a maximal system of islands and related problems are presented.
In the present paper, we list related problems with exact formulae.
In each case, we present the proof, which we believe to be the shortest.

In full generality, we denote the set of all cells of some board by $\mathcal C$.
A {\it height function} is a  mapping $h: \mathcal C\rm\to \mathbb R$, $c \mapsto h(c)$.
We have to specify a neighborhood relation on the cells.
If not otherwise stated, two cells are {\it neighbors} if they share a point.
Let $R$ be a subset of cells.
The neighbors of $R$ can be defined naturally as the set of cells not in $R$ but
having a neighbor in $R$.
A connected subset $R$ of cells is called an {\it island}, if
the minimum height in $R$ is greater than the maximum height on the neighbors of $R$.
In our applications, we define the islands to have a geometric shape,
therefore the definition of connectivity does not play a role here.
If $h$ is a height function, then we denote the induced set of islands  by ${\mathcal I}(h)$.
Let us consider rectangular islands.
We say that rectangles $R$ and $S$ are \it far from each other, \rm if no cell of $R$ is the neighbor of any cell of $S$.
We denote by $P(\mathcal C)$ the power set of $\mathcal C$, that is the set of all subsets of $\mathcal C$.
The following statement in a different form was proved in \cite{Cz}.

\begin{lemma} \label{egy}
Let $\mathcal C$ be the set of all cells of some board, and let $\mathcal W$ denote the entire board as an island.
Let ${\mathcal I}$ be a set of islands.
The following two conditions are equivalent:
\item{\rm(i)} there exists a mapping $h: \mathcal C\rm\to \mathbb R$, $c \mapsto h(c)$ such that
${\mathcal I}  = {\mathcal I}(h)$.
\item{\rm (ii)} $\mathcal B \rm \in {\mathcal I}$, and for any
$R_1\neq R_2 \in {\mathcal I}$ either $R_1\subset R_2$, or $R_2\subset R_1$,
or $R_1$ and $R_2$ are far from each other.
\end{lemma}

A subset of $P(\mathcal C)$ satisfying the equivalent conditions
of Lemma~\ref{egy} is called \it a system of islands\rm .
The set of maximal elements of
 $\mathcal I\setminus \{\mathcal B\rm \}$
is denoted by $\max \mathcal I$.


\section{Methods}\label{meth}

We list three effective proof techniques for island problems.
We give detailed demonstration of the latter two, the first and original method can be read in \cite{Cz}.

We recall the following

\begin{lemma}[\cite{Cz}]
The maximum number of rectangular islands of an $m \times n $  rectangular board is
$$f(m,n)= \left[\frac{(m+1)(n+1)}{2}\right]-1.$$
\end{lemma}

Let $\mathcal C$ be the set of unit squares of the $m \times n $ board.
The proof in \cite{Cz} exploits that the islands form a weakly independent set
in the distributive lattice of $P(\mathcal C)$.
In a distributive lattice, maximal weakly independent subsets are called \it weak bases\rm .
By the main theorem of \cite{CzHSch}, any two weak bases have the same cardinality.
We ask the reader to consult \cite{Cz} for the details.

For the second method, we need basic graph theory \cite{bondy}.
To be self-contained, we recall the definitions that are crucial for our purposes.
A graph without a cycle is called a {\it forest}.
Any component of a forest is a connected cycle-free graph, that is a {\it tree}.
A forest with a distinguished node (root) in each component is called a {\it rooted forest}.
For any node $u$, there is a unique path from $u$ to the root of its component.
If $u$ is not a root, then this path has more than one vertex.
Let $u^+$ be the node following $u$ on the path to the root.
It is called the {\it father of $u$}.
If $v=u^+$, then we say that $u$ is a son of $v$.
If $v$ is on the path from $u$ to a root, then $v$ is an {\it ancestor} of
$u$ and $u$ is a {\it descendant} of $v$.
Any non-root vertex has exactly one father, but a father might have several sons.
The descendants of $v$ are the sons of $v$, the sons of sons (grandsons),
and so on.
For any $v$ the vertex $v$ and its descendants span $T_v$, a {\it
rooted subtree}.
Therefore, a rooted forest can be described recursively:
it contains a set of roots; each root has a set of sons; 
and there are vertex disjoint rooted trees rooted at the sons.
A vertex is a {\it leaf} if and only if it has no son.
A rooted tree is {\it binary}
if and only if any non-leaf node has two sons.

Consider any base set.
In the present paper, it is the set $\mathcal C$ of all cells.
Fix certain shapes (e.g. rectangle) to be allowed for islands, and a function $h$ defined on $\mathcal C$.
Let $\mathcal I$ be the set of islands of the fixed shape.

\begin{fact}
Let $S$ be a subset of $\mathcal C$.
The maximal islands contained in $S$ are disjoint.
\end{fact}

Based on this observation, we define a rooted forest, $T({\mathcal I})$ describing a
hierarchy of the islands.
Let the maximal islands $R_1, R_2,\ldots, R_t$ of $\mathcal I$ be the roots of the forest.
The islands contained in $R$ form $P(R)$
($R\in P(R)$), the part of the partition connected to $R$.
The maximal islands of $P(R)-\{R\}$ are the sons of $R$.
The description of the rooted forest is completed by iterating the above step.

\begin{remark}
In the specific cases we consider,
the base set is always an island itself, therefore it is the unique maximal island.
In this case, the rooted forest is a rooted tree.
\end{remark}

Let $T_0({\mathcal I})$ be the rooted forest, we just defined
based on $\mathcal I$. The islands are exactly the vertices of
$T_0({\mathcal I})$, hence the number of islands is
$|V(T_0({\mathcal I})|$. The leaves of $T_0({\mathcal I})$
are the minimal islands.

We can visualize this description.
The function $h$ can be viewed as a height function, describing
a geographic part of Earth. We start to pour water into this place.
We see the birth of a few islands, these are the roots.
Often water level zero corresponds to the case when we see the first
island: the whole area we considered.
As the water level increases, we see islands to be divided into smaller islands (sons) or disappear (leaves of our forest).

Sometimes an island/vertex has only one son.
This means that by the increase of the water level the island
shrinks.
In this case, it will be useful to modify our rooted forest.
We interpret the decline of the island as a division
into a smaller island (its only son) and a dummy part.
This dummy part of the island will
be a second son of the shrinking vertex, a leaf.
Let $T({\mathcal I})$ be the rooted forest we obtain this way.
In $T({\mathcal I})$ any non-leaf vertex has at least two sons.
The number of islands is $|V(T({\mathcal I}))|-|D|$,
where $D$ is the set of dummy nodes added to $T_0({\mathcal I})$.
The leaves of $T({\mathcal I})$ are the minimal islands and the dummy islands.

\begin{figure}[ht]
 \begin{center}
  \includegraphics[scale=0.5]{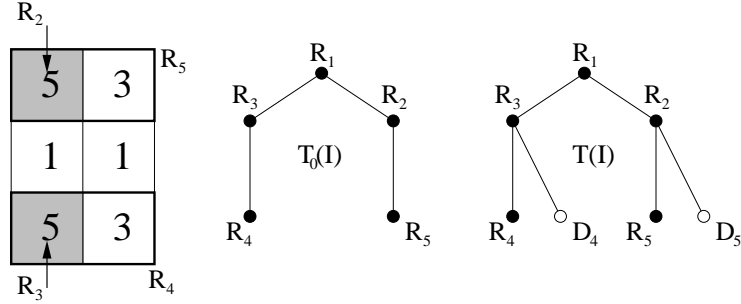}
 \end{center}
   \caption{Hasse diagram of islands with respect to containment}
  \label{tree}
\end{figure}

We demonstrate the above notation in Figure~\ref{tree}.

In order to bound the number of islands, the following Lemma (folklore or an easy exercise in studying rooted trees) is very useful.

\begin{lemma} \label{ketto}
\item{(i)} Let $T$ be a binary tree with $\ell$ leaves. Then
the number of vertices of $T$ depends only on $\ell$ and
$|V|=2\ell-1$.
\item{(ii)} Let $T$ be a rooted tree such that any non-leaf node has
at least two sons. Let $\ell$ be the number of leaves in $T$.
Then $|V|\leq 2\ell -1$.
\end{lemma}

Our simple strategy is the following:
if we know how to express the number of islands by the number of vertices and
dummy nodes, then we apply Lemma~\ref{ketto}.

\begin{proof}[A proof example]
Let $B_{m,n}$ denote the set of $mn$ unit squares of the $m\times n$ rectangular board.
Let the island shape be rectangular. 
We call the vertices of the unit squares grid points.

Let $\mathcal I$ be a system of islands with $s$ minimal islands and $d$
dummy islands.
Any island covers at least four grid points.
In the case of a shrinking island there is a loss of at least two grid points.
The set of these lost grid points can be assigned to the corresponding dummy node.
We assign grid points to the leaves of $T({\mathcal I})$: 
four points to the minimal islands, two points to the dummy leaves.
These assigned sets of  grid points are disjoint in the set of all $(m+1)(n+1)$ grid points.
Therefore, $4s+2d\leq (m+1)(n+1)$.
The number of leaves of $T({\mathcal I})$ is $\ell=s+d$.
By  Lemma~\ref{ketto}, the number of islands
is $|V|-d\leq (2\ell-1)-d=2s+d-1\leq \frac{1}{2}(m+1)(n+1)-1$.
\end{proof}

This proof is very suggestive, clear and short.
Still, it needed some technical preparation.
As it turns out, we can make the proof even more elementary.

The iterative description of $T_0({\mathcal I})$ or
$T({\mathcal I})$ suggests a recursive proof technique:
the mathematical induction.
Actually, all known upper bounds on the number of islands \cite{Cz, HNP,  P} can be proved by induction.

\begin{proof}[A proof example]
Let $f(B_{m,n})=f(m,n)$ be the maximum number of islands on the $m\times n$ rectangular board.
We claim that $f(B_{m,n})\leq \frac{1}{2}(m+1)(n+1)-1$.
Let us denote the covered grid points by $\|B_{m,n}\|$.
For disjoint sub-boards $S_1, S_2, \ldots S_k$ of $B_{m,n}$
we know that $\|B_{m,n}\|\geq \|S_1\|+\|S_2\|+\ldots +\|S_k\|$ holds.

We prove the claim by induction.
The case of small boards can be easily checked.
Let $\mathcal I^*$ be a system of islands
realizing the number $f(m,n)$.

\begin{multline*}
f(m,n)= 1+\sum_{R\in{\max \mathcal I^*}} f(R)\leq 1+\sum_{R\in{\max \mathcal I^*}}\left(\frac{1}{2}\|R\|-1\right)=\\
= 1+\frac{1}{2}\sum_{R\in{\max \mathcal I^*}}\|R\|-|{\max \mathcal I^*}| \leq \frac{1}{2}\|B_{m,n}\|+1-|{\max \mathcal I^*}|.
\end{multline*}

If $|{\max \mathcal I^*}|\geq 2$, then the induction is complete.
If $|{\max \mathcal I^*}|=1$, then one needs a minor technical remark to finish the proof.
\end{proof}


\section{Applications} \label{app}
\subsection{Peninsulas}

We show that the maximum cardinality of rectangular islands in the
 $ m \times n $  rectangular board can be attained such that 
 each island reaches at least one side of the board. 
 This is a slight strengthening of the result in \cite{Cz}.
 Also, the proof gives  a recursive algorithm constructing a system of maximum cardinality.
For brevity, we call a rectangular island $P$ a {\it peninsula}
if it reaches at least one side of the board. 
We denote the maximum number of peninsulas in an $m \times n$ board by $p(m,n)$.

\begin{thm}\label{pmn}
In a rectangular $m\times n$ board, the maximum number of rectangular islands is equal to the maximum number of peninsulas, that is $p(m,n)=f(m,n)$.
\end{thm}

{\it Proof.}
Since peninsulas are islands,  $p(m,n)\leq f(m,n)$.
To prove $p(m,n)\geq f(m,n)$,
we show by induction on the number of cells, that the maximum number of peninsulas reaching the eastern 
side of the board is at least $f(m,n)$.
We use the notation $p'(m,n)$ for the  maximum number of peninsulas reaching the
eastern side of the board.
For $m,n\in \{1,2\}$, the statement is clear.
To see the induction step, notice the following:
Let the first row of the board be a peninsula.
It contains $m$ different peninsulas by deleting the squares one by one from west.
That is,

$$p'(m,n)\geq p'(m,n-2)+m=\left[\frac{m(n-2)+m+n-2-1}{2}\right] +m+1= \left[\frac{mn+m+n-1}{2}\right].\qed $$


\subsection{Cylindric board, rectangular islands}\label{cyli1}

In this section, we put a square grid on the surface of a cylinder with height $m$ and circumference of the base circle $n$.
We get the same object by identifying the sides of length $m$ of an $m\times n$ rectangle.
We denote by $c_1(m,n)$ the maximum number of rectangular islands on this cylinder, supposing that the whole cylinder is an
island, but no other cylinders are islands.

\begin{thm}\label{c1mn}
If $n\geq 2$, then $c_1(m,n)=\left[\frac{(m+1)n}{2}\right].$
\end{thm}

{\it Proof.} By deleting a column of the cylinder, we get an $m\times (n-1)$ rectangle. Therefore,

$$c_1(m,n)\geq f(m, n-1)+1=\left[\frac{(m+1)n}{2}\right].$$

Let $\mathcal I^*$ be a  set of rectangular islands of maximum cardinality.
Then

\begin{multline*}
c_1(m,n)= 1+\sum_{R\in{max \mathcal I^*}} f(R)= 1+\sum_{R\in{max\mathcal I^*}}\left( \left[\frac{(u+1)(v+1)}{2}  \right]-1 \right)=\\
=1 - | max ({\mathcal I^*})| +\sum_{R\in max \mathcal I^*} \left[\frac{(u+1)(v+1)}{2}  \right]\le\\
\le 1 - 1 + \left[\frac{(m+1)n}{2} \right] = \left[\frac{(m+1)n}{2} \right] .
\end{multline*}

We applied that $-| max ({\mathcal I^*})|$ can be bounded above
by $-1$ if $| max ({\mathcal I^*})| \geq 1$; and also that

$$\sum_{R\in max \mathcal I^*} \left[\frac{(u+1)(v+1)}{2}  \right]\le \left[\frac{(m+1)n}{2} \right].$$

\begin{figure}[ht]
 \begin{center}
\includegraphics[scale=0.4]{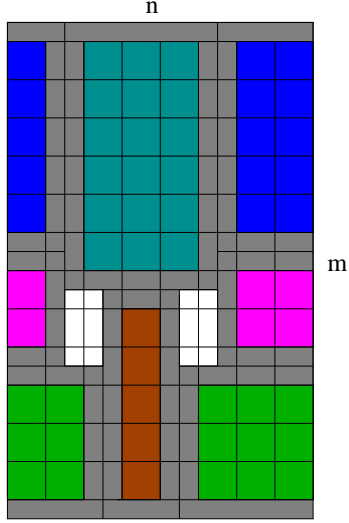}
   \caption{Magnified rectangles of a cylindric board}
  \label{hengernagyit}
 \end{center}
\end{figure}

To see this, we magnify the maximal rectangles by half a unit, see
Figure~\ref{hengernagyit}.
The sum of the area of the magnified maximal rectangles
is at most the area of the magnified cylinder. \qed

\subsection{Cylindric board, cylindric and rectangular islands}

Living on a cylindric board, it is natural to consider cylindric islands as well.
In this section, we allow two shapes for the islands, cylindric and rectangular.
We denote by $c_2(m,n)$ the maximum cardinality of such a system of islands on the cylindric $m\times n$ board.

 \begin{thm}\label{c2mn}
If $n\geq 2$, then $c_2(m,n)=\left[\frac{(m+1)n}{2}\right]+\left[\frac{(m-1)}{2}\right].$
\end{thm}

{\it Proof.}
We show by induction on $m$, that $c_2(m,n)\geq \left[\frac{(m+1)n}{2}\right]+\left[\frac{m-1}{2}\right]$.
Notice that $c_2(1,n)=n$ and $c_2(2,n)\geq f(2,n-1)+1=\left[\frac{3n}{2}\right].$
Let $m>2$. 
For the induction step, we remove a cylinder of height one such that a 
cylindric board of size $(m-2)\times n$ remains.
Therefore,

$$c_2(m,n)\geq c_2(m-2,n)+n+1= \left[\frac{(m-1)n}{2}\right]+\left[\frac{m-3}{2}\right]+n+1=\left[\frac{(m+1)n}{2}\right]+\left[\frac{m-1}{2}\right].$$

Now we show that $c_2(m,n)\leq \left[\frac{(m+1)n}{2}\right]+\left[\frac{m-1}{2}\right]$.
There must be a cylindric island $Y$ by Theorem~\ref{c1mn}.
Then $Y$ is included in a maximal cylindric island $M$. 
Assume there is a maximum cardinality system given.
Then $M$ is bordered with water from one side, and in the maximal case the width of this water is one.
Also, the remaining part of the board is a maximal cylindric island. 
Therefore, there exist $a,b\in \mathbb N_0$ in such a way that
$a+b+1=m$ and

\begin{multline*}  c_2(m,n)= c_2(a,n)+c_2(b,n)+1 = \left[\frac{(a+1)n}{2}\right]+\left[\frac{a-1}{2}\right]+ \left[\frac{(b+1)n}{2}\right]+\left[\frac{b-1}{2}\right]+ 1  \leq \\
\leq \left[\frac{(a+b+1+1)n}{2}\right]+\left[\frac{a+b-3+1}{2}\right] +1
=\left[\frac{(m+1)n}{2}\right]+\left[\frac{m-1}{2}\right]. \qed
\end{multline*}


\subsection{Toroidal board, rectangular islands}

With respect to the neighborhood relation, the most symmetric case is the toroidal board.
It is not a surprise that we get the most compact result of all.

Assume there is an $m \times n$ board on the torus, which is also known as $C_m\times C_n$.
The island shape is fixed as rectangular, but we consider the whole torus as an island.
We denote by $t(m,n)$ the maximum number of rectangular islands on the torus.

\begin{thm}\label{tmn}
 If $m, n\geq 2$, then $t(m,n)=\left[\frac{mn}{2}\right].$
\end{thm}

{\it Proof.}
We can cut off a horizontal and a vertical line to get an $(m-1)\times (n-1)$ rectangle.
Therefore,

$$t(m,n)\geq f(m-1,n-1)+1=\left[\frac{mn}{2}\right]. $$

On the other hand,
we denote again by $\mathcal I^*$ a set of rectangular  islands realizing the maximum cardinality.

\begin{multline*}  t(m,n)= 1+\sum_{R\in max \mathcal I^*} f(R) = 1+\sum_{R\in max \mathcal I^*} \left( \left[\frac{(u+1)(v+1)}{2}  \right]-1 \right) =\\
= 1 - | max ({\mathcal I^*})| +\sum_{R\in \max \mathcal I^*} \left[\frac{(u+1)(v+1)}{2}  \right]
\leq   1 - 1 + \left[\frac{mn}{2}\right] = \left[\frac{mn}{2}\right].
\end{multline*}

Similar to Section~\ref{cyli1} we applied that  $-| max ({\mathcal I^*})|$ can be bounded above
by $-1$ if $| max ({\mathcal I^*})| \geq 1$; and also that

$$\sum_{R\in \max \mathcal I^*} \left[\frac{(u+1)(v+1)}{2}  \right]\leq \left[\frac{mn}{2}\right], $$ 

by counting the grid points covered by maximal islands. \qed

\subsection{Heuristic}

The results of the section are based on counting grid points.
Therefore, it is convenient to modify the parameters in the above cases such that
the number of grid points of the board is the same.
This yields: 
$$p(m-1,n-1)=f(m-1,n-1)=c_1(m-1,n)-1=t(m,n)-1.$$

We can erase the $-1$ if the sets of grid points induce the islands 
instead of the squares of the board.

We could short-cut to this result as follows:
The most restricted set of islands is the eastern peninsulas and the
broadest is the toroidal definition. 
Therefore, $p(m-1,n-1)\le f(m-1,n-1)\le c_1(m-1,n)\le t(m,n)$. 
Then we observe that a peninsula construction coincides with the
toroidal upper bound. That is, equality must hold everywhere.


\section{Changing the neighborhood relation}

So far, two cells were neighbors if they had a point in common.
Therefore, in the corresponding neighborhood graph, the typical degree was 8.
It is somewhat natural to rarefy this structure such that the
neighborhood graph is also a grid.
In this case, two cells are neighbors if and only if they have a side in common.

\begin{figure}[ht]
 \begin{center}
\includegraphics[scale=0.4]{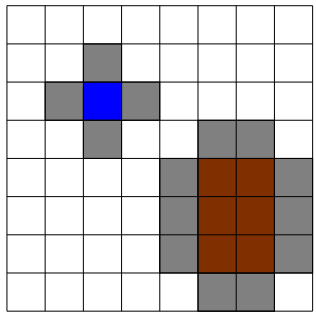}
  \caption{Common-side neighborhood}
  \label{szomi}
 \end{center}
\end{figure}

Let a rectangular board of size $m \times n$ be given.
We denote the maximum cardinality of a system of rectangular islands by $\hat{f} (m,n)$.
We denote the set of rectangular islands induced by a height function $h$  by ${\hat{\mathcal I}}(h)$.

\begin{lemma}
Let $\mathcal C$ be the set of cells of an $m \times n$ board, and let $\mathcal W$
denote the entire board as an island.
Let ${\mathcal I}$ be a subset of\/ $P(\mathcal C)$.
Two cells are neighbors if and only if they have a side in common.
The following two conditions are equivalent:
\item{\rm(i)} there exists a mapping $h: \mathcal C\rm\to \mathbb R$, $c \mapsto h(c)$ such that
${\mathcal I}  = {\hat{\mathcal I}}(h)$.
\item{\rm (ii)} $\mathcal W \rm \in {\mathcal I}$, and for any
$I_1\neq I_2 \in {\mathcal I}$ either $I_1\subset I_2$, or $I_2\subset I_1$,
or $I_1$ and $I_2$ have zero or one point in common.
\end{lemma}

Despite the rarefied structure, the maximum value surprisingly has not changed.

\begin{thm}\label{hatf}
 $\hat{f} (m,n)=f(m,n)$.
\end{thm}

\begin{figure}[ht]
 \begin{center}
\includegraphics[scale=0.4]{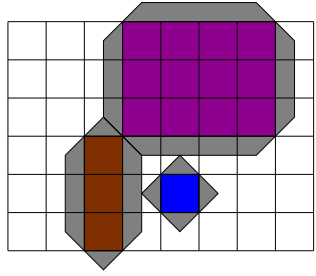}
  \caption{Magnification of the maximal islands}
  \label{szu}
 \end{center}
\end{figure}

{\it Proof.}
Clearly $\hat{f} (m,n)\geq f(m,n)$.
We show $\hat{f} (m,n)\leq f(m,n)$ via induction on $mn$.
If $m,n\in \{1,2\}$, then the statement holds.
We denote by\ $\mathcal I^*$ a maximum cardinality system of islands.
For the induction step, we magnify the members of $max \mathcal I^*$ as
shown in Figure~\ref{szu}.

Let the side lengths of a rectangle $R$ be $u$ and $v$.
We define $\mu(R)=\mu(u,v):=(u+1)(v+1)-2$, which is the area of the magnified $R$.
The induction step goes as follows

\begin{multline*}
\hat{f}(m,n) =   1+\sum_{R\in max \mathcal I^*} \hat{f}(R)= 1+\sum_{R\in max \mathcal I^*} \left( \left[\frac{(u+1)(v+1)}{2}\right]-1 \right)  = \\
=   1 + \sum_{R\in max \mathcal I^*} \left( \left[\frac{\mu(u,v)}{2}\right]\right)
  \leq 1 + \left[ \frac{\mu(C)}{2} \right] -1.
\end{multline*}

In the last step we applied that
$$\sum_{R\in max \mathcal I^*} {\mu(u,v)}  \leq \mu (\mathcal C) -2, $$
since the magnified maximal islands do not overlap, and
there is an area of size at least 2, which is not covered by the magnified maximal islands.
\qed

\section{Islands in hypercubes}

We give an exact formula for the maximum number of hypercubic islands in a big hypercube.
The board consists of all vertices of a  hypercube, or in other words the elements of a Boolean algebra $\mathcal A=\{0,1\}^n$.
Two cells are neighbors if their Hamming distance is 1.
We denote the maximum number of islands in $\mathcal A=\{0,1\}^n$ by $b(n)$.

\begin{thm}\label{bn}
$b(n)= 1 + 2^{n-1}.$
\end{thm}

{\it Proof.}
Consider the vertices with an odd number of 1's.
They correspond to independent cellular islands.
Therefore, $b(n)\geq 1+2^{n-1}$, if we consider the whole space as an island.

We prove the opposite direction by induction on $n$.
For $n=0,1$ the statement is easy to check.
For $n\geq 2$, we cut the hypercube into two half-hypercubes of size $2^{n-1}$.
If one of them is an island, then the other part can not contain an island.
If neither of them is an island, then by the induction hypothesis, in both half-hypercubes, the maximum cardinality of a system of islands is at most $2^{n-2}$.
This implies the claim:
$b(n)\leq 1+2^{n-1}$. \qed


\section*{Epilogue}

The inductive argument worked easily, when we found a row or column containing no island.
This phenomenon helped us also when we traced back the toroidal case to the planar.
Let an empty row or column be called a {\it blast}.
The maximum cardinality of a blast-free system of islands can be of interest.
(Blast-free domino tilings of a rectangular board is a classical Olympiad question.)
As this maximum is strongly related to the area uncovered by the maximal islands,
it is tempting to ask the following

\begin{prob}
Let us consider the $m\times n$ rectangular, cylindric or toroidal board.
What is the minimum of the uncovered area in a blast-free configuration of
maximal islands?
\end{prob}

We dare to conjecture $m+n+1$ in the plane,
$3m+2n-7$ on the cylinder and $4m+2n-9$ on the torus, where $m\le n$.

We may assume the maximal islands to be on the same level, height one say.
This is the first level also in the rooted tree defined in Section~\ref{meth}.
We imagine the islands corresponding to vertices of the second
level of the rooted tree to have height two.
This building process can be continued downwards the rooted tree.
In this way, we build a characteristic example of a class of island systems corresponding
to the same rooted tree.
This can also be formulated in the language of the height function.

In this sense, the previous question is posed about the section of a landscape at height one.
We may require the blast-free property at each level.

\begin{prob}
Let $H=\max_{c\in \mathcal C} h(c)$.
For which height function $h$ are the sections blast-free at each height $1,2,\dots, H$. 
\end{prob}

Another natural question is the following

\begin{prob}
Characterize the maximum cardinality island systems on any fixed board.
\end{prob}

The philosophy of Section~\ref{app} can be applied to grid-like drawings of 
orientable surfaces of higher genus as well as non-orientable surfaces.

We plan to report on these problems in a forthcoming paper.


\end{document}